# Limit theorems for functionals on the facets of stationary random tessellations

LOTHAR HEINRICH[1], HENDRIK SCHMIDT[2] and VOLKER SCHMIDT[3]

[1]*Institute of Mathematics, University of Augsburg, 86135 Augsburg, Germany.*
*E-mail: lothar.heinrich@math.uni-augsburg.de*
[2]*France Telecom NSM/R&D/RESA/NET, 92794 Issy Moulineaux, Cedex 9, France.*
*E-mail: hendrik.schmidt@orange-ftgroup.com*
[3]*Institute of Stochastics, Ulm University, 89069 Ulm, Germany.*
*E-mail: volker.schmidt@uni-ulm.de*

We observe stationary random tessellations $X = \{\Xi_n\}_{n \geq 1}$ in $\mathbb{R}^d$ through a convex sampling window $W$ that expands unboundedly and we determine the total $(k-1)$-volume of those $(k-1)$-dimensional manifold processes which are induced on the $k$-facets of $X$ ($1 \leq k \leq d-1$) by their intersections with the $(d-1)$-facets of independent and identically distributed motion-invariant tessellations $X_n$ generated within each cell $\Xi_n$ of $X$. The cases of $X$ being either a Poisson hyperplane tessellation or a random tessellation with weak dependences are treated separately. In both cases, however, we obtain that all of the total volumes measured in $W$ are approximately normally distributed when $W$ is sufficiently large. Structural formulae for mean values and asymptotic variances are derived and explicit numerical values are given for planar Poisson–Voronoi tessellations (PVTs) and Poisson line tessellations (PLTs).

## 1. Introduction

In this paper, we consider stationary random tessellations $X = \{\Xi_n\}_{n \geq 1}$ of the $d$-dimensional Euclidean space $\mathbb{R}^d$ with convex cells $\Xi_n$. We assume that within each cell $\Xi_n$ of the initial tessellation $X$, a further random tessellation $X_n = \{\Xi_{n\ell}\}_{\ell \geq 1}$ of $\mathbb{R}^d$ is nested, that is, $\Xi_n$ is subdivided into cells $\Xi_n \cap \Xi_{n\ell}, \ell \geq 1$, where the sequence of component tessellations $(X_n)_{n \geq 1}$ consists of independent copies of a generic motion-invariant tessellation $X_0$ drawn independently of $X$. The assumption of motion-invariance of $X_0$ will play a crucial role in deriving explicit moment formulae. This type of iterated random tessellation is said to be an $X/X_0$-*nesting* in $\mathbb{R}^d$. Having available only a single observation of such an $X/X_0$-nesting in a presumably large, convex sampling window $W$, we are interested in the asymptotic behaviour of the random sums

$$Z_k^{(d)}(W) = \sum_{n \geq 1} \vartheta_n^{(k)}(W), \qquad 1 \leq k \leq d-1, \tag{1.1}$$

where 'asymptotic' means that $W \uparrow \mathbb{R}^d$ and where the random measures

$$\vartheta_n^{(k)}(\cdot) = \nu_{k-1}(X_n^{(d-1)} \cap \Xi_n^{(k)} \cap \cdot), \qquad 1 \leq k \leq d-1, \tag{1.2}$$

act on the Borel sets of $\mathbb{R}^d$. The functional $\vartheta_n^{(k)}(W)$ measures the $(k-1)$-volume of the random subsets induced in $W$ by the intersection of the motion-invariant manifold process $X_n^{(d-1)}$ of $(d-1)$-facets of $X_n$ with the union $\Xi_n^{(k)}$ of all $k$-faces belonging to the boundary $\partial\Xi_n$ of the $n$th cell of $X$ (cf. Section 2.1 for details and precise definitions). Note that, by definition, the random measures $\vartheta_1^{(k)}(\cdot), \vartheta_2^{(k)}(\cdot), \ldots$ in (1.2) are conditionally independent given the tessellation $X$.

Our results supplement earlier central limits theorems (CLTs) for cumulative measures of stationary ergodic tessellations modelling the total effect of random internal cell structures ([14]). Whereas, in the latter reference, the random measures corresponding to those in the sum (1.1) act on the interiors of the cells $\Xi_n$, the measures $\vartheta_n^{(k)}(\cdot)$ defined in (1.2) are concentrated on the cell boundaries $\partial\Xi_n$ of $X$. Hence, certain new effects arise due to the interactions between the stationary random manifold process $\bigcup_{n \geq 1} \partial\Xi_n$ of cell boundaries of $X$ and the component tessellations $(X_n)_{n \geq 1}$. It turns out that there are considerable differences between $X$ being a stationary Poisson hyperplane tessellation (PHT) and $X$ satisfying certain weak dependence assumptions. In the first case, due to the overnormalization in the CLTs for Poisson hyperplane processes caused by inherent long-range dependences, (cf. [15]), the influence of $X_0$ on the Gaussian limit distribution is relatively weak. The other case seems to be somewhat more delicate because the asymptotic variance of the existing Gaussian limits are influenced by first and second order characteristics of both $X$ and $X_0$.

We present our derived CLTs in the general case of $\mathbb{R}^d$, since this allows for a clearer and more transparent exposition. Clearly, however, the CLTs find their applications in the modelling of planar, but also spatial, networks as they occur, for example, in cell biology and telecommunications. Indeed, concentrating on the latter example, the problem often arises of handling and modelling data that represent the geometrical structure of the infrastructure system (e.g., main roads and side streets) that supports the technical telecommunications equipment. In recent years, stochastic-geometric modelling approaches have proven useful and are established domains of research today. In particular, the Stochastic Subscriber Line Model (SSLM) has been developed as an integrated and easily extendable model for telecommunication access networks (cf. [25] and the references therein).

The SSLM employs (iterated) random tessellations to describe the geometric network support. Having identified the best fitting model from a class of potentially suitable tessellations (cf. [8]), cost functionals and their distributions can be studied along the network geometry (cf. [9]).

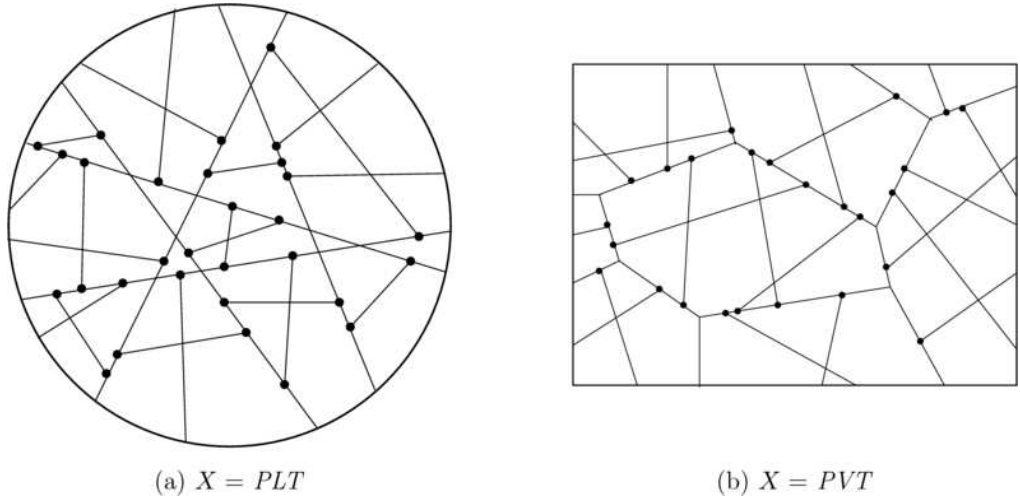

(a) $X = PLT$                                        (b) $X = PVT$

**Figure 1.** Realizations of planar tessellations $X$, where the nested tessellation $X_0$ is a PLT. The T-crossings are displayed as thick dots ●.

Assume that we use a planar $X/X_0$-nesting to model the geometric support. In the framework of our study, we observe a single value $Z_1^{(2)}(W)$, which counts the number of T-crossings in a sampling region $W \subset \mathbb{R}^2$ induced by the intersection of the edges of the tessellation $X_n$ with the edges of the $n$th cell $\Xi_n$ for $n \geq 1$. Figure 1 shows two examples of this situation in differently shaped sampling windows $W$. In particular, Figure 1(a) shows a PLT/PLT-nesting through a ball of radius $r > 0$ and centered at the origin, whereas, in Figure 1(b), we consider a PVT/PLT-nesting within a rectangular sampling window.

The analysis of the aforementioned T-crossings, that is, of the connections between main roads and side streets, plays an important role in telecommunication modelling since these crossings are the entry points to the blockwise civil engineering of the local network. Let the type of the initial tessellation $X$ and the type of the nested tessellation $X_0$ be known. Within a suitably large region $W$, the distribution of the number of T-crossings is then known through our results. Thus, the engineer is provided with useful information about the local network. For example, it is possible to deduce the dimensioning and capacity potential for each entry point in order to provide blockwise optimal connection quality to the subscribers, where one block comprises all of those subscribers who are situated in the cells formed by the main roads.

In contrast to that, assume that we have, again in a suitably large region $W$, knowledge about local information, such as the type of $X_0$, and especially about the value $Z_1^{(2)}(W)$ for the T-crossings. The expression on the left-hand side of (4.1), as well as the expression on the left-hand side of (5.7), can then be calculated and used to test for normality.

Depending on the (unknown) type of $X$ (representing the main road system), we expect to reject the null hypothesis of normality either for the formula in (4.1) or for the formula in (5.7). This can provide, in the framework of model selection, a hint as to the structure of $X$ before passing to more refined fitting procedures.

The paper is organized as follows. In Section 2, we introduce basic notation and recall some relevant facts from stochastic geometry. Section 3 presents mean value relations and formulae for (asymptotic) variances. In Sections 4 and 5, we formulate and prove the announced CLTs for the different cases of initial tessellations $X$. Finally, in Section 6, we study some examples of weakly dependent tessellations and discuss possible extensions of our results.

## 2. Preliminaries

In this section, we introduce the basic notation and present a brief account of some relevant material on random tessellations and stochastic geometry in general. For a detailed and rigorous discussion of these topics, we refer to the existing mathematical literature, in particular to [21, 22, 24, 27] and [30], which contain many further references, as well as numerous tessellation models with applications to various fields.

Throughout, let $(\Omega, \sigma(\Omega), \mathbb{P})$ be a common probability space on which all random objects occurring in the present paper will be defined. Further, let $\langle x, y \rangle = \sum_{k=1}^{d} x_k y_k$ denote the scalar product of the coordinate vectors $x = (x_1, \ldots, x_d)^\top$ and $y = (y_1, \ldots, y_d)^\top$ in $\mathbb{R}^d$. By means of the Euclidean norm $\| \cdot \| = \sqrt{\langle \cdot, \cdot \rangle}$, we define the closed ball $B_r^d = \{x \in \mathbb{R}^d : \|x\| \le r\}$ with radius $r \ge 0$ centered at the origin and the unit sphere $\mathbb{S}^{d-1} = \{x \in \mathbb{R}^d : \|x\| = 1\}$ in $\mathbb{R}^d$, respectively. Remember that each affine $(d-1)$-dimensional subspace $H$ of $\mathbb{R}^d$, called *hyperplane* in $\mathbb{R}^d$ in the sequel, admits a parameter representation $H(p, v) = \{x \in \mathbb{R}^d : \langle x, v \rangle = p\}$. Here, $p \in \mathbb{R}^1$ denotes the signed perpendicular distance of $H$ from the origin and $v \in \mathbb{S}_+^{d-1} = \{(x_1, \ldots, x_d)^\top \in \mathbb{S}^{d-1} : x_d \ge 0\}$ is the directional vector belonging to the upper unit hemisphere. Further, let $\nu_k(\cdot)$ denote the Lebesgue or $k$-volume measure in $\mathbb{R}^k$ for $k = 0, \ldots, d$, where we can also just write $\nu_d(\cdot) = | \cdot |$. The $k$-dimensional Lebesgue measure will also be used instead of the $k$-dimensional Hausdorff measure on (affine) $k$-dimensional subspaces in $\mathbb{R}^d$ for any $k = 0, \ldots, d-1$. As usual, $\nu_0(\cdot)$ coincides with the counting measure, that is, $\nu_0(B) = \#B$. For brevity, put

$$\kappa_d = |B_1^d| = \frac{\pi^{d/2}}{\Gamma(d/2 + 1)}, \qquad \text{where } \Gamma(s) = \int_0^\infty \mathrm{e}^{-y} y^{s-1}\, \mathrm{d}y \text{ for } s > 0.$$

The families of all non-empty closed, compact and compact convex sets in $\mathbb{R}^d$ are denoted by $\mathcal{F}_d'$, $\mathcal{K}_d'$ and $\mathcal{C}_d'$, respectively. Note that $\mathbb{B}(S)$ stands for the $\sigma$-algebra of Borel sets in the metric space $S$.

## 2.1. Random tessellations and random nestings

In this section, we sketch out the mathematically rigorous approach to random tessellations, as used in stochastic geometry, and we recall some basic facts, where referring [21, 22, 24, 27] and [30] for a systematic study of these topics.

A tessellation of $\mathbb{R}^d$ is a countable family $\tau = \{C_n\}_{n \geq 1}$ of convex bodies $C_n \in \mathcal{C}'_d$ such that $\operatorname{int} C_n \neq \emptyset$ for all $n$, $\operatorname{int} C_n \cap \operatorname{int} C_m = \emptyset$ for all $n \neq m$, $\bigcup_{n \geq 1} C_n = \mathbb{R}^d$ and $\sum_{n \geq 1} \mathbb{1}_{\{C_n \cap K \neq \emptyset\}} < \infty$ for any $K \in \mathcal{K}'_d$. Notice that the sets $C_n$, called the *cells* of $\tau$, are necessarily polytopes in $\mathbb{R}^d$. The family of all tessellations in $\mathbb{R}^d$ is denoted by $\mathcal{T}$. A random tessellation $X = \{\Xi_n\}_{n \geq 1}$ in $\mathbb{R}^d$ is a sequence of random convex bodies $\Xi_n$ such that $\mathbb{P}(X \in \mathcal{T}) = 1$.

Note that a (stationary) random tessellation $X$ can also be modelled as a (stationary) marked point process $\sum_{n \geq 1} \delta_{[\alpha(\Xi_n), \Xi_n^0]}$, where $\alpha \colon \mathcal{C}'_d \to \mathbb{R}^d$ is a $\mathcal{B}(\mathcal{F}'_d)$-measurable mapping such that $\alpha(C) \in C$ and $\alpha(C + x) = \alpha(C) + x$ for any $C \in \mathcal{C}'_d$ and $x \in \mathbb{R}^d$, and where $\Xi_n^0 = \Xi_n - \alpha(\Xi_n)$ is the centered cell corresponding to $\Xi_n$ which contains the origin. The point $\alpha(C)$ is called the *associated point* of $C$ and is usually chosen to be the centroid or lexicographically smallest point of $C$.

Suppose that the stationary marked point process $\sum_{n \geq 1} \delta_{[\alpha(\Xi_n), \Xi_n^0]}$ has positive and finite intensity $\gamma = \mathrm{E}\#\{n \colon \alpha(\Xi_n) \in [0, 1)^d\}$. By $\mathcal{P}_d^0$, we denote the set of all compact and convex $d$-polytopes whose associated point is located at the origin. The corresponding Palm mark distribution $P^0$ of $X$ is then given by

$$P^0(B) = \gamma^{-1} \mathrm{E}\#\{n \colon \alpha(\Xi_n) \in [0, 1)^d, \Xi_n^0 \in B\}, \qquad B \in \mathcal{B}(\mathcal{F}'_d) \cap \mathcal{P}_d^0. \tag{2.1}$$

The notion of a *typical cell* of $X$ refers to a random polytope $\Xi^* \colon \Omega \to \mathcal{P}_d^0$ whose distribution coincides with $P^0$. Since the cells $\Xi_n$ are space-filling and non-overlapping (up to a null set), we have the mean value relationship

$$\frac{1}{\gamma} = \int_{\mathcal{P}_d^0} |C| P^0(\mathrm{d}C), \tag{2.2}$$

that is, the cell intensity $\gamma$ equals the reciprocal of $\mathrm{E}|\Xi^*|$.

A (deterministic) iterated tessellation $\tau = \{C_{n\ell} \cap C_n \colon \operatorname{int} C_{n\ell} \cap \operatorname{int} C_n \neq \emptyset\}$ in $\mathbb{R}^d$ consists of an initial tessellation $\tau = \{C_n\}_{n \geq 1}$ in $\mathbb{R}^d$ and a sequence $(\tau_n)_{n \geq 1}$ of component tessellations $\tau_n = \{C_{n\ell}\}_{\ell \geq 1}$. In order to define a random iterated tessellation, we proceed along the lines of [19]. Let $\Xi$ be a random convex body in $\mathbb{R}^d$ with $\mathbb{P}$-a.s. non-empty interior and let $X = \{\Xi_n\}_{n \geq 1}$ be a random tessellation in $\mathbb{R}^d$. Then, the counting measure $Y(\cdot \mid \Xi)$ defined by $Y(B \mid \Xi) = \sum_{n \geq 1} \delta_{\Xi_n \cap \Xi}(B) \mathbb{1}_{\{\operatorname{int} \Xi_n \cap \operatorname{int} \Xi \neq \emptyset\}}$ for $B \in \mathcal{B}(\mathcal{F}'_d)$ is a point process in $\mathcal{C}'_d$ describing a random tessellation of $\Xi$.

Furthermore, if $X = \{\Xi_n\}_{n \geq 1}$ is an arbitrary random tessellation in $\mathbb{R}^d$ and if $X_n = \{\Xi_{n\ell}\}_{\ell \geq 1}$, $n = 1, 2 \dots$, are independent copies of a generic random tessellation $X_0$ in $\mathbb{R}^d$ drawn independent of $X$, then the random counting measure $Y(B) = \sum_n Y_n(B \mid \Xi_n)$, where $Y_n(B \mid \Xi_n) = \sum_{\ell \geq 1} \delta_{\Xi_{n\ell} \cap \Xi_n}(B) \mathbb{1}_{\{\operatorname{int} \Xi_{n\ell} \cap \operatorname{int} \Xi_n \neq \emptyset\}}$ for $B \in \mathcal{B}(\mathcal{F}'_d)$, is called the *point-process representation* of an iterated random tessellation (briefly , an $X/X_0$-*nesting*) in $\mathbb{R}^d$ with initial tessellation $X$ and component tessellations $X_1, X_2, \ldots$. Clearly, the point process $Y$ is stationary (and isotropic), provided that both the initial tessellation $X$ and the generic component tessellation $X_0$ are stationary (and isotropic). Moreover, $Y$ is ergodic if $X$ possesses this property.

Each stationary (motion-invariant) random tessellation $X = \{\Xi_n\}_{n \geq 1}$ in $\mathbb{R}^d$ induces $d$ stationary (motion-invariant) random lower-dimensional manifold processes $X^{(k)}$, called *$k$-facet processes* of $X$ for $k = 0, 1, \ldots, d - 1$. For example, $X^{(0)}$ is the point process of *vertices* and $X^{(1)}$ is the line segment process of *edges* of $X$.

To be precise, $X^{(k)}$ is defined to be the union of all of the *$k$-facets* of $X$, whereas $\Xi_n^{(k)}$ denotes the union of all *$k$-faces* of its $n$th cell $\Xi_n$. Here, the $k$-facets of $X$ are $k$-polytopes in $\mathbb{R}^k$ which arise from a finite intersection of neighbouring cells of $X$. The $(d-1)$-faces of $\Xi_n$ are $(d-1)$-polytopes in the boundary $\partial \Xi_n$ and $k$-faces are defined recursively for $k = 0, \ldots, d - 2$ as $k$-polytopes in the relative boundaries of the $(k+1)$-faces. Note that the set of all $k$-faces may differ from the set of $k$-facets and that, for example in [27], Chapter 6, $X^{(k)}$ is used slightly differently to denote the point process of $k$-facets.

A random tessellation $X = \{\Xi_n\}_{n \geq 1}$ in $\mathbb{R}^d$ is said to be *normal* (or *ordinary*) if $\mathbb{P}$-a.s. every $k$-facet of $X$ lies in the boundaries of exactly $d - k + 1$ cells, $k = 0, \ldots, d - 1$. Many-real life tessellations in $\mathbb{R}^2$ and $\mathbb{R}^3$ possess this property, which motivates the term 'normal'. There are important classes of stationary tessellations in $\mathbb{R}^d$ whose cells are constructed (realizationwise), according to specific geometric rules, from the atoms of a stationary point process in $\mathbb{R}^d$. Among them are *Voronoi* and *Laguerre* tessellations (see, e.g., [24] for details), which turn out to be normal if the generating point process is Poisson; see [21]. It seems that this fact continues to hold for a large class of (even instationary) generating point processes which are mixing in a certain sense and/or whose higher-order moment measures possess Lebesgue densities. In [13], it is shown that Voronoi tessellations in $\mathbb{R}^d$ are normal if the $(d+2)$th-order product density of the generating stationary point process exists. Computations of second-order characteristics of spatial Poisson–Voronoi tessellations can be found in [11]. Finally, it should be mentioned that there are more general definitions of tessellations (cf., e.g., [32]), allowing for the rigorous treatment of random tessellations which do not necessarily consist of only convex cells. Without doubt, the most prominent example is the Johnson–Mehl tessellation; see also [24] for details. For this model, CLTs have been proved, based on $\alpha$-mixing conditions derived from the generating (Poisson) point process (cf. [4] and [5]).

## 2.2. Stationary Poisson hyperplane tessellations

Let $\Psi = \sum_{i \geq 1} \delta_{[P_i, V_i]}$ be a stationary and independently marked Poisson point process on the real line $\mathbb{R}^1$ with intensity $\lambda$ and mark distribution $\Theta$ on the mark space $\mathbb{S}_+^{d-1}$; see [6]. By means of the parameter representation $H(p, v)$, $(p, v) \in \mathbb{R}^1 \times \mathbb{S}_+^{d-1}$, of a hyperplane in $\mathbb{R}^d$, we may represent a Poisson hyperplane process (PHT) $\Phi$ (defined in [27] as a point process on the space of affine $(d-1)$-dimensional subspaces in $\mathbb{R}^d$) with intensity $\lambda$ and

(spherical) orientation distribution $\Theta$ by

$$\Phi = \sum_{i \geq 1} \delta_{H(P_i, V_i)}. \tag{2.3}$$

The Poisson hyperplane process $\Phi$ given in (2.3) is said to be *non-degenerate* if $\Theta(H(0,v) \cap \mathbb{S}_+^{d-1}) < 1$ for any $v \in \mathbb{S}_+^{d-1}$. In this case, (2.3) induces stationary *$k$-flat processes* $\Phi_k$ for $k = 0, 1, \ldots, d-1$ whose countable support consists of the affine $k$-dimensional subspaces ($k$- intersection flats) $H(P_{i_1}, V_{i_1}) \cap \cdots \cap H(P_{i_{d-k}}, V_{i_{d-k}})$ for pairwise distinct indices $i_1, \ldots, i_{d-k} \geq 1$. The union of these $k$-flats coincides with the $k$-facet process $X^{(k)}$ of the corresponding stationary PHT $X = \{\Xi_n\}_{n \geq 1}$ generated by (2.3). The cells $\Xi_n$, $n \geq 1$, are bounded $d$-polytopes ($\mathbb{P}$-a.s.) if and only if $\Phi$ is non-degenerate; see [27], Chapter 6. Furthermore, this property implies that the stationary $k$-volume measure $\vartheta_{k,d}(\cdot)$ associated with $\Phi_k$ (resp. $X^{(k)}$) and defined by

$$\vartheta_{k,d}(B) = \frac{1}{(d-k)!} \sum_{i_1, \ldots, i_{d-k} \geq 1}^{*} \nu_k \left( \bigcap_{j=1}^{d-k} H(P_{i_j}, V_{i_j}) \cap B \right) \qquad \text{for bounded } B \in \mathcal{B}(\mathbb{R}^d), \tag{2.4}$$

where $\sum^{*}$ denotes summation over pairwise distinct indices, has positive intensity

$$\lambda_{k,d} = \mathbb{E}\vartheta_{k,d}([0,1]^d) = \frac{(2\lambda)^{d-k}}{(d-k)!\kappa_d} \mathbb{E}g_{k,d}(Q_0, V_0) \qquad \text{for } k = 0, 1, \ldots, d-1. \tag{2.5}$$

Here, the function $(p,v) \mapsto g_k^{(d)}(p,v)$ is defined on $[-1,1] \times \mathbb{S}_+^{d-1}$ by

$$g_{k,d}(p,v) = \mathbb{E}\nu_k \left( \bigcap_{i=1}^{d-k-1} H(Q_i, V_i) \cap H(p,v) \cap B_1^d \right), \tag{2.6}$$

where $(Q_i, V_i)$, $i = 0, 1, \ldots, d-1$, are i.i.d. random vectors with independent components. Note that the generic random variable $Q_0$ is uniformly distributed on $[-1,1]$ and the generic random vector $V_0$ has the orientation distribution $\Theta$; see [15].

It is well known from convex geometry (see [26]) Chapter 3.5, that the probability measure $\Theta$ on $\mathbb{S}_+^{d-1}$ determines a unique, centrally symmetric convex body $Z_\Theta$, called the *associated zonoid*, which is given by

$$h(Z_\Theta, u) = \int_{\mathbb{S}_+^{d-1}} |\langle u, v \rangle| \Theta(dv) \qquad \text{for } u \in \mathbb{R}^d,$$

where $h(K, u) = \max_{x \in K} \langle u, x \rangle$ denotes the support function of an arbitrary $K \in \mathcal{C}_d'$.

In [16], the following closed-form expression of $g_{k,d}(p,v)$ in terms of $Z_\Theta$ has been derived

$$g_{k,d}(p,v) = \frac{(d-k-1)!\kappa_{d-1}}{2^{d-k-1}} (1-p^2)^{(d-1)/2} \mathbb{1}_{[-1,1]}(p) V_{d-k-1}^{(d-1)}(Z_\Theta^v), \tag{2.7}$$

where $K^v$ denotes the image of $K \in \mathcal{C}'_d$ under orthogonal projection onto $H(0, v)$ and $V_j^{(d-1)}(K)$ stands for the intrinsic $j$-volume of $K \in \mathcal{C}'_{d-1}$. Using the relationship

$$jV_j^{(d)}(Z_\Theta) = \int_{\mathbb{S}_+^{d-1}} V_{j-1}^{(d-1)}(Z_\Theta^v)\Theta(\mathrm{d}v) \qquad \text{for } j = 1, \ldots, d$$

(cf. [26, 31], Chapter 3.5 and [16]) combined with (2.5) and (2.7) yields that $\lambda_{k,d} = \lambda^{d-k} \times V_{d-k}^{(d)}(Z_\Theta)$ for $k = 0, 1, \ldots, d-1$, which has already been stated in [20], Chapter 6.

The stationary Poisson hyperplane process $\Phi$ given in (2.3) is isotropic (and hence motion-invariant) if and only if $\Theta$ is the uniform distribution, which means that $Z_\Theta = \frac{\kappa_{d-1}}{d\kappa_d}B_1^d$. This, in turn, leads to the explicit formula

$$\lambda_{k,d} = \binom{d}{k}\frac{\kappa_d}{\kappa_k}\left(\frac{\kappa_{d-1}}{d\kappa_d}\right)^{d-k}\lambda^{d-k} \qquad \text{for } k = 0, 1, \ldots, d-1. \tag{2.8}$$

We are now in a position to formulate a CLT for the total $k$-volume $\vartheta_{k,d}(B_\varrho^d)$ of the support of the $k$-flat process $\Phi_k$ contained in the ball $B_\varrho^d$. This result has been proven in [15], even in a multidimensional version.

**Theorem 2.1.** *Let* $\Phi = \sum_{i \geq 1} \delta_{H(P_i, V_i)}$ *be a stationary, non-degenerate Poisson hyperplane process with orientation distribution* $\Theta$ *on* $\mathbb{S}_+^{d-1}$ *and intensity* $\lambda > 0$*. Then,*

$$\frac{\vartheta_{k,d}(B_\varrho^d) - \lambda_{k,d}|B_\varrho^d|}{|B_\varrho^d|^{1-1/2d}} \xrightarrow[\varrho \to \infty]{d} \mathcal{N}(0, \sigma_{k,d}^2) \qquad \text{for } k = 0, 1, \ldots, d-1, \tag{2.9}$$

*where*

$$\sigma_{k,d}^2 = \lim_{\varrho \to \infty} \frac{\mathrm{Var}(\vartheta_k^{(d)}(B_\varrho^d))}{|B_\varrho^d|^{2-1/d}} = \frac{(2\lambda)^{2d-2k-1}}{((d-k-1)!)^2\kappa_d^{2-1/d}}\mathbb{E}g_{k,d}^2(Q_0, V_0), \tag{2.10}$$

*with* $g_{k,d}(p, v)$ *and* $(Q_0, V_0)$ *defined by* (2.6)*. If, additionally,* $\Phi$ *is isotropic, that is,* $\Theta$ *is the uniform distribution on* $\mathbb{S}_+^{d-1}$*, then* $\lambda_{k,d}$ *is given by* (2.8) *and* $\sigma_{k,d}^2$ *takes the explicit form*

$$\sigma_{k,d}^2 = \lambda^{2d-2k-1}\frac{2^{2d-1}\kappa_d^{1/d}}{(2d-1)!}\binom{d-1}{k}^2\left(\frac{d!\kappa_{d-1}}{k!\kappa_k}\right)^2\left(\frac{\kappa_{d-1}}{d\kappa_d}\right)^{2(d-k)}. \tag{2.11}$$

Note that, even in the anisotropic case, we have

$$\sigma_{d-1,d}^2 = \lambda\frac{2^{2d-1}\kappa_{d-1}^2}{(2d-1)!\kappa_d^{2-1/d}}, \tag{2.12}$$

that is, $\sigma_{d-1,d}^2$ coincides with the left-hand side of (2.11) for $k = d-1$. This is due to the fact that the $(d-1)$-volume of the hyperplanes $H(P_i, V_i)$ within the ball $B_\varrho^d$ does not depend on $V_i$ and so the distribution of $\vartheta_{k,d}(B_\varrho^d)$ is independent of the orientation distribution $\Theta$. Furthermore, we mention that in [16], Theorem 2.1 could be extended to non-spherical convex sampling windows $W_\varrho = \varrho W_1$ (cf. Section 5 below). However, in this case, the formulae (2.7) and (2.11) depend on $W_1$ and are less explicit.

## 3. First- and second-order moment formulae

Let $X = \{\Xi_n\}_{n \geq 1}$ be a stationary random tessellation of $\mathbb{R}^d$ and let $X_0$ be a motion-invariant tessellation independent of $X$. We consider an $X/X_0$-nesting in $\mathbb{R}^d$, as in Section 2.1, observed within a convex sampling window $W$. In order to calculate expectation and variance of the random variables $Z_k^{(d)}$ in (1.1) for $k = 1, \ldots, d-1$, we need two intensity values.

First, we consider $\lambda_0^{(k,d)}$, the intensity of the stationary $(k-1)$-dimensional manifold process $X_0^{(d-1)} \cap L$ generated by the intersection of the $(d-1)$-facet process $X_0^{(d-1)}$ with an arbitrary $k$-flat $L$ in $\mathbb{R}^d$. Since $X_0^{(d-1)}$ is motion-invariant by the assumption of motion-invariance of $X_0$, we may identify $L$ with $\mathbb{R}^k$ so that $\lambda_0^{(k,d)}$ can be defined by

$$\lambda_0^{(k,d)} = \mathbb{E}\nu_{k-1}(X_0^{(d-1)} \cap \mathbb{R}^k \cap [0,1)^k), \qquad 1 \leq k \leq d-1. \qquad (3.1)$$

By using quite general stereological relationships derived in [18], we may express $\lambda_0^{(k,d)}$ by the (full-dimensional) intensity $\lambda_0^{(d,d)} = \mathbb{E}\nu_{d-1}(X_0^{(d-1)} \cap [0,1)^d)$ of the manifold process of $(d-1)$-facets $X_0^{(d-1)}$ through

$$\lambda_0^{(k,d)} = c_k^{(d)} \lambda_0^{(d,d)} \qquad \text{with } c_k^{(d)} = \frac{\gamma((k+1)/2)\gamma(d/2)}{\gamma(k/2)\gamma((d+1)/2)}, \ 1 \leq k \leq d-1. \qquad (3.2)$$

Further, let

$$\mu_k^{(d)} = \mathbb{E}\nu_k(X^{(k)} \cap [0,1)^d) \qquad (3.3)$$

denote the intensity of the stationary $k$-facet process $X^{(k)}$ associated with $X$. To avoid rather involved formulae, in particular, for the variance of $Z_k^{(d)}(W)$, we impose an additional condition on the tessellation $X = \{\Xi_n\}_{n \geq 1}$.

**Condition F.** *For $k = 1, \ldots, d-1$, assume that there exists a non-random integer $m_k^{(d)} \geq 1$ such that*

$$m_k^{(d)} \nu_k(X^{(k)} \cap W) = \sum_{n \geq 1} \nu_k(\Xi_n^{(k)} \cap W) \qquad \mathbb{P}\text{-}a.s.$$

*for any $W \in \mathcal{C}_d'$ with $|W| > 0$.*

Condition F means that, for $k = 1, \ldots, d-1$, each $k$-facet of $X$ lies in a constant number $m_k^{(d)}$ of $k$-faces of cells $\Xi_n, n \geq 1$. Obviously, Condition F is satisfied for any planar tessellation with $m_1^{(2)} = 2$. For $d \geq 3$, by the very definition, any normal tessellation $X$ obeys Condition F with $m_k^{(d)} = d - k + 1$, and a non-degenerate stationary Poisson hyperplane tessellation $X$ is easily seen to satisfy Condition F with $m_k^{(d)} = 2^{d-k}$, see [27]. Note, however, that Poisson–Delaunay tessellations do not satisfy Condition F for $d \geq 3$.

**Lemma 3.1.** *Consider an $X/X_0$-nesting in $\mathbb{R}^d$ with stationary initial tessellation $X = \{\Xi_n\}_{n \geq 1}$ and motion–invariant component tessellation $X_0$. Assume that $X$ satisfies Condition F and that $0 < \gamma^{-1} = \mathbb{E}|\Xi^*| < \infty$ (cf. (2.2)). If $\mu_k^{(d)} < \infty$ or, equivalently, $\mathbb{E}\nu_k(\Xi^{*(k)}) < \infty$, and $\lambda_0^{(k,d)} < \infty$ for any $k = 1, \ldots, d-1$, then*

$$\mathbb{E}Z_k^{(d)}(W) = \lambda_0^{(k,d)} m_k^{(d)} \mu_k^{(d)} |W| \qquad \text{for any } W \in \mathcal{C}_d', \quad k = 1, \ldots, d-1. \tag{3.4}$$

*Moreover, if, additionally,*

$$\mathbb{E}\nu_k^2(X^{(k)} \cap [0,1)^d) < \infty \quad \text{and} \quad \int_{\mathcal{P}_d^0} \mathbb{E}\nu_{k-1}^2(X_0^{(d-1)} \cap C^{(k)})P^0(\mathrm{d}C) < \infty, \tag{3.5}$$

*then for $W \in \mathcal{C}_d'$ and $k = 1, \ldots, d-1$,*

$$\begin{aligned}
\mathrm{Var}(Z_k^{(d)}(W)) = \gamma \int_{\mathcal{P}_d^0} \int_{\mathbb{R}^d} \mathrm{Var}(\nu_{k-1}(X_0^{d-1} \cap C^{(k)} \cap (W - x))) \, \mathrm{d}x P^0(\mathrm{d}C) \\
+ (\lambda_0^{(k,d)} m_k^{(d)})^2 \, \mathrm{Var}(\nu_k(X^{(k)} \cap W)).
\end{aligned} \tag{3.6}$$

**Proof.** Let $k \in \{1, \ldots, d-1\}$ be fixed and let $\mathbb{E}_X(\cdot)$ denote the conditional expectation $\mathbb{E}(\cdot|X)$ given the tessellation $X = \{\Xi_n\}_{n \geq 1}$. Hence, we may rewrite the expectation of $Z_k^{(d)}(W)$ introduced in (1.1) as

$$\mathbb{E}Z_k^{(d)}(W) = \mathbb{E}\sum_{n \geq 1} \mathbb{E}_X \nu_{k-1}(X_0^{(d-1)} \cap \Xi_n^{(k)} \cap W).$$

Owing to the motion-invariance of $X_0^{(d-1)}$, we get, together with (3.1), that

$$\mathbb{E}_X \nu_{k-1}(X_0^{(d-1)} \cap \Xi_n^{(k)} \cap W) = \lambda_0^{(k,d)} \nu_k(\Xi_n^{(k)} \cap W)$$

for any cell $\Xi_n$. In view of Condition F, we may proceed by writing that

$$\mathbb{E}_X Z_k^{(d)}(W) = \lambda_0^{(k,d)} \sum_{n \geq 1} \nu_k(\Xi_n^{(k)} \cap W) = \lambda_0^{(k,d)} m_k^{(d)} \nu_k(X^{(k)} \cap W). \tag{3.7}$$

Combined with the stationarity of $X^{(k)}$, this gives $\mathbb{E}\nu_k(X^{(k)} \cap W) = \mu_k^{(d)}|W|$, which, in turn, proves (3.4). Recall that by using the notion of the typical cell (cf. (2.2)), we have

$$\mathbb{E}\sum_{n \geq 1} \nu_k(\Xi_n^{(k)} \cap W) = \gamma \mathbb{E}\nu_k(\Xi^{*(k)})|W|(= m_k^{(d)} \lambda_{k,d}|W|),$$

which establishes the relationship $\mathbb{E}\nu_k(\Xi^{*,(k)}) = m_k^{(d)}\lambda_{k,d}\mathbb{E}|\Xi^*|$. To verify (3.6), we start with the well-known identity

$$\mathrm{Var}(Z_k^{(d)}(W)) = \mathbb{E}(\mathrm{Var}_X\, Z_k^{(d)}(W)) + \mathrm{Var}(\mathbb{E}_X Z_k^{(d)}(W)), \tag{3.8}$$

where $\mathrm{Var}_X(\cdot)$ denotes the conditional variance $\mathrm{Var}(\cdot|X)$ given $X$. Since, conditional on the tessellation $X = \{\Xi_n\}_{n\geq 1}$ the random measures $\vartheta_1^{(k)}(\cdot), \vartheta_2^{(k)}(\cdot), \ldots$ in (1.2) are stochastically independent, we obtain

$$\mathrm{Var}_X\, Z_k^{(d)}(W) = \sum_{n\geq 1} \mathrm{Var}_X(\nu_{k-1}(X_0^{(d-1)} \cap \Xi_n^{(k)} \cap W)). \tag{3.9}$$

With $\Xi_n = \Xi_n^0 + \alpha(\Xi_n)$, we may apply the refined Campbell theorem to the stationary marked point process $\sum_{n\geq 1} \delta_{[\alpha(\Xi_n), \Xi_n^0]}$ (cf. [6] or [21]), where we find, together with (2.1), that

$$\mathbb{E}(\mathrm{Var}_X\, Z_k^{(d)}(W)) = \gamma \int_{\mathbb{R}^d} \int_{\mathcal{P}_d^0} \mathrm{Var}(\nu_{k-1}(X_0^{(d-1)} \cap (C^{(k)} + x) \cap W)) P^0(\mathrm{d}C)\, \mathrm{d}x$$

$$= \gamma \int_{\mathcal{P}_d^0} \int_{\mathbb{R}^d} \mathrm{Var}(\nu_{k-1}(X_0^{(d-1)} \cap C^{(k)} \cap (W - x)))\, \mathrm{d}x\, P^0(\mathrm{d}C). \tag{3.10}$$

Here, we used (3.5), the stationarity of $X_0^{(d-1)}$, Fubini's theorem and the fact that

$$\nu_{k-1}((B + x) \cap W) = \nu_{k-1}(B \cap (W - x))$$

for any bounded $B \in \mathcal{B}(\mathbb{R}^{k-1})$ and $x \in \mathbb{R}^d$. The existence of the inner Lebesgue integral in the second line of (3.10) is also seen by applying Fubini's theorem and the second condition of (3.5), that is,

$$\int_{\mathbb{R}^d} \mathbb{E}\nu_{k-1}^2(X_0^{(d-1)} \cap C^{(k)} \cap (W - x))\, \mathrm{d}x$$

$$= \mathbb{E} \int_{X_0^{(d-1)} \cap C^{(k)}} \int_{X_0^{(d-1)} \cap C^{(k)}} |(W - u) \cap (W - v)|\nu_{k-1}(\mathrm{d}u)\nu_{k-1}(\mathrm{d}v)$$

$$\leq \mathbb{E}\nu_{k-1}^2(X_0^{(d-1)} \cap C^{(k)})|W|.$$

From (3.7), combined with the first condition in (3.5), it is immediately seen that the second term on the right-hand side of (3.8) is finite and takes the form

$$\mathrm{Var}(\mathbb{E}_X Z_k^{(d)}(W)) = (\lambda_0^{(k,d)} m_k^{(d)})^2\, \mathrm{Var}(\nu_k(X^{(k)} \cap W)).$$

The latter equality, together with (3.10) and (3.8), confirms the validity of (3.6).    $\square$

The second condition of (3.5) imposes restrictions on both the initial and the component tessellation. Note that this condition is fulfilled if

$$\mathbb{E}\nu_{k-1}^2(X_0^{(d-1)} \cap \mathbb{R}^k \cap [0,1)^k) < \infty \quad \text{and} \quad \mathbb{E}N_k^2(\Xi^*)(1+D(\Xi^*))^{2k} < \infty \quad (3.11)$$

for $k = 1, \ldots, d-1$, where $N_k(C)$ and $D(C) = \sup\{\|x-y\| : x, y \in C\}$ denote the number of $k$-faces and the diameter of the $d$-polytope $C \in \mathcal{P}_d^0$, respectively. To see that (3.11) implies the second condition in (3.5), we write $C^{(k)}$ as union of the $k$-faces $C_l^{(k)}$, $l = 1, \ldots, N_k(C)$, and use the motion invariance of $X_0^{(d-1)}$ to obtain the estimate

$$\begin{aligned}
\mathbb{E}\nu_{k-1}^2(X_0^{(d-1)} \cap C^{(k)}) &\leq \mathbb{E}\left(\sum_{l=1}^{N_k(C)} \nu_{k-1}(X_0^{(d-1)} \cap C_l^{(k)})\right)^2 \\
&\leq N_k^2(C)\mathbb{E}\nu_{k-1}^2(X_0^{(d-1)} \cap \mathbb{R}^k \cap [0,D(C))^k) \\
&\leq N_k^2(C)(1+D(C))^{2k}\mathbb{E}\nu_{k-1}^2(X_0^{(d-1)} \cap \mathbb{R}^k \cap [0,1)^k)
\end{aligned}$$

for $k = 1, \ldots, d-1$.

Note that if $\vartheta_n^{(k)}(\cdot)$ acts on the interior of $\Xi_n$, as supposed in [14], then the conditional expectation $\mathbb{E}_X \vartheta_n^{(k)}(W)$ is a constant multiple of $|\Xi_n \cap W|$ and therefore $\mathbb{E}_X Z_k^{(d)}(W) = \sum_{n\geq 1} \mathbb{E}_X \vartheta_n^{(k)}(W)$ is proportional to $|W|$. In this case, the second variance term on the right-hand side of (3.8) vanishes. Due to this fact, the formula for the variance of $Z_k^{(d)}(W)$ given in [14] is relatively simple and the ergodicity of the initial tessellation $X$ suffices to prove asymptotic normality of $Z_k^{(d)}(W)$.

# 4. CLTs for manifold processes on facets of a stationary PHT

In this section, we consider the random measures $\vartheta_n^{(k)}(\cdot)$ given in (1.2) for $n \geq 1$ whose support lies in the cell boundaries $\partial\Xi_n$, more precisely in the $k$-facets (for $k = 1, \ldots, d-1$) of a non-degenerate stationary PHT $X = \{\Xi_n\}_{n\geq 1}$ in $\mathbb{R}^d$. For the sake of simplicity, we assume in this section that the sampling window $W$ is the $d$-dimensional ball $B_\varrho^d$ centered at the origin and with radius $\varrho > 0$. For more general expanding sampling windows, we refer to the comment at the end of Section 2.2. By Theorem 2.1, combined with the geometric and probabilistic properties of PHTs, we prove a CLT (as $\varrho \to \infty$) for the total $(k-1)$-volume of the sets (contained in $B_\varrho^d$) arising from the intersection of the $k$-faces of $\Xi_n$ with the $(d-1)$-facets of the component tessellations $X_n$ for $n \geq 1$.

**Theorem 4.1.** *Let $X = \{\Xi_n\}_{n\geq 1}$ be the stationary PHT in $\mathbb{R}^d$ generated by a stationary non-degenerate Poisson hyperplane process $\Phi$ given in (2.3) with orientation distribution $\Theta$ and intensity $\lambda > 0$. Furthermore, let $X_0$ be a motion-invariant random tessellation in $\mathbb{R}^d$ having the intensities $\lambda_0^{(k,d)} > 0$ (cf. (3.1) and (3.2), resp.) Assume* that $\mathbb{E}\nu_{k-1}^2(X_0^{(d-1)} \cap \mathbb{R}^k \cap (0,1]^k) < \infty$ for $k = 1, \ldots, d-1$ and that the corresponding $X/X_0$-nesting is observed through the spherical sampling window $B_\varrho^d$. Then,

$$\frac{Z_k^{(d)}(B_\varrho^d) - \widetilde{\mu}_k^{(d)}|B_\varrho^d|}{|B_\varrho^d|^{1-1/(2d)}} \underset{n\to\infty}{\longrightarrow} \mathcal{N}(0, \widetilde{\sigma}_{k,d}^2) \qquad \text{for } k = 1, \ldots, d-1, \tag{4.1}$$

where

$$\widetilde{\mu}_k^{(d)} = 2^{d-k}\lambda_0^{(k,d)}\lambda_{k,d}, \qquad \widetilde{\sigma}_{k,d}^2 = \lim_{\varrho\to\infty}\frac{\mathrm{Var}(Z_k^{(d)}(B_\varrho^d))}{|B_\varrho^d|^{2-1/d}} = (2^{d-k}\lambda_0^{(k,d)})^2\sigma_{k,d}^2 \tag{4.2}$$

and $\lambda_{k,d}$, $\sigma_{k,d}^2$, and $\lambda_0^{(k,d)}$ are defined by (2.5), (2.10) and (3.1), respectively.

**Proof.** We first recall that in case of a stationary PHT $X$, we have that $m_k^{(d)} = 2^{d-k}$ and that the intensity (3.3) of the $k$-facet process $X^{(k)}$ coincides with the intensity (2.5) of the $k$-flat process $\Phi_k$ induced by (2.3), that is, we have that $\mu_k^{(d)} = \lambda_{k,d}$. Hence, the formulae for the intensities $\widetilde{\mu}_k^{(d)} = \mathbb{E}Z_k^{(d)}([0,1)^d)$ of $Z_k^{(d)}(\cdot)$ follow from (3.4), as stated in (4.2). Next, we rewrite the mean zero random variable $Z_k^{(d)}(B_\varrho^d) - \widetilde{\mu}_k^{(d)}|B_\varrho^d|$ as

$$Z_k^{(d)}(B_\varrho^d) - \widetilde{\mu}_k^{(d)}|B_\varrho^d| = S_\varrho^{(k)} + 2^{d-k}\lambda_0^{(k,d)}T_\varrho^{(k)}, \tag{4.3}$$

where

$$S_\varrho^{(k)} = Z_k^{(d)}(B_\varrho^d) - \mathbb{E}_X Z_k^{(d)}(B_\varrho^d) = Z_k^{(d)}(B_\varrho^d) - 2^{d-k}\lambda_0^{(k,d)}\nu_k(X^{(k)} \cap B_\varrho^d)$$

and

$$T_\varrho^{(k)} = \nu_k(X^{(k)} \cap B_\varrho^d) - \mu_k^{(d)}|B_\varrho^d| = \vartheta_{k,d}(B_\varrho^d) - \lambda_{k,d}|B_\varrho^d|.$$

From Theorem 2.1, we obtain that

$$\frac{T_\varrho^{(k)}}{|B_\varrho^d|^{1-1/(2d)}} \underset{n\to\infty}{\longrightarrow} \mathcal{N}(0, \sigma_{k,d}^2) \qquad \text{for } k = 1, \ldots, d-1. \tag{4.4}$$

By means of Slutsky's theorem (cf., e.g., [17]), the proof of the CLT (4.1) is complete whenever

$$\frac{S_\varrho^{(k)}}{|B_\varrho^d|^{1-1/(2d)}} \underset{\varrho\to\infty}{\overset{P}{\longrightarrow}} 0. \tag{4.5}$$

In view of Chebychev's inequality, we need only to prove that

$$\frac{\mathbb{E}(S_\varrho^{(k)})^2}{|B_\varrho^d|^{2-1/d}} \underset{\varrho\to\infty}{\longrightarrow} 0, \tag{4.6}$$

which includes first of all to ensure that $\mathbb{E}(S_\varrho^{(k)})^2 < \infty$. Since $\mathbb{E}(S_\varrho^{(k)})^2 = \mathbb{E}\operatorname{Var}_X(Z_\varrho^{(k)})$, we obtain, in analogy to the proof of Lemma 3.1 and by taking into account (3.11), that

$$\mathbb{E}(S_\varrho^{(k)})^2 = \gamma \int_{\mathcal{P}_d^0} \int_{\mathbb{R}^d} \operatorname{Var}(\nu_{k-1}(X_0^{d-1} \cap C^{(k)} \cap (B_\varrho^d - x)))\,dx\,P^0(dC)$$

$$\leq \gamma \int_{\mathcal{P}_d^0} \mathbb{E}\nu_{k-1}^2(X_0^{(d-1)} \cap C^{(k)})P^0(dC)|B_\varrho^d|$$

$$\leq \gamma \mathbb{E}\nu_{k-1}^2(X_0^{(d-1)} \cap \mathbb{R}^k \cap [0,1)^k)\mathbb{E}N_k^2(\Xi^*)(1 + D(\Xi^*))^{2k}|B_\varrho^d|.$$

Using the distributional properties of the typical cell $\Xi^*$ of a stationary PHT, in particular that $D(\Xi^*)$ has an exponential moment (cf. [2]), we find that

$$\mathbb{E}N_k^2(\Xi^*)(1 + D(\Xi^*))^{2k} < \infty \qquad \text{for } k = 1, \ldots, d-1,$$

which immediately confirms (4.6) for any $d \geq 2$ (cf. also [7] and [28]). Finally, using the formula (3.6) for the variance of $Z_k^{(d)}(B_\varrho^d)$, together with the limiting relations (4.6) and (2.10), we find that

$$\widetilde{\sigma}_{k,d}^2 = \lim_{\varrho \to \infty} \frac{\operatorname{Var}(Z_k^{(d)}(B_\varrho^d))}{|B_\varrho^d|^{2-1/d}} = (2^{d-k}\lambda_0^{(k,d)})^2\sigma_{k,d}^2.$$

This completes the proof of Theorem 4.1. $\qquad\qquad\qquad\qquad\qquad\qquad\square$

# 5. CLTs for manifold processes induced on the facets of a stationary weakly dependent tessellation

Throughout this section, we consider a stationary $X/X_0$-nesting which can be observed through an expanding family of convex sampling windows $W_\varrho$ with shape $W_\varrho = \varrho W_1$ for $\varrho > 0$, where $W_1 \in \mathcal{C}_d'$ contains a ball and is itself contained in a ball, that is, $B_r^d \subseteq W_1 \subseteq B_R^d$ for some $0 < r \leq R < \infty$. We assume that the stationary initial tessellation $X = \{\Xi_n\}_{n \geq 1}$ is ergodic (cf. [6, 24, 27]) and possesses, in contrast to Poisson hyperplane tessellations, further weak dependence properties. The latter properties ensure asymptotic normality of the total $k$-volume of the $k$-facets in a large sampling window $W_\varrho$. More precisely, we impose on $X$ the following condition.

**Condition G.** *For $k = 1, \ldots, d-1$, assume that there exists a real number $\tau_{k,d}^2 \geq 0$ such that*

$$\frac{\operatorname{Var}(\nu_k(X^{(k)} \cap W_\varrho))}{|W_\varrho|} \xrightarrow[\varrho \to \infty]{} \tau_{k,d}^2$$

*and*

$$\frac{\nu_k(X^{(k)} \cap W_\varrho) - \mu_k^{(d)}|W_\varrho|}{|W_\varrho|^{1/2}} \xrightarrow[n \to \infty]{\mathrm{d}} \mathcal{N}(0, \tau_{k,d}^2).$$

In analogy to Section 4, we shall prove that the centered and normalized cumulative functional (1.1) on $W_\varrho$, that is,

$$\frac{Z_k^{(d)}(W_\varrho) - \eta_k^{(d)}|W_\varrho|}{|W_\varrho|^{1/2}} \qquad \text{with } \eta_k^{(d)} = \mathbb{E}Z_k^{(d)}([0,1)^d), \tag{5.1}$$

converges in distribution to a Gaussian random variable $\mathcal{N}(0, \widetilde{\tau}_{k,d}^2)$ for $k = 1, \ldots, d-1$, as $W_\varrho \uparrow \mathbb{R}^d$. Here, $\widetilde{\tau}_{k,d}^2$ denotes the asymptotic variance of the random variable (5.1) as $\varrho \to \infty$, that is

$$\widetilde{\tau}_{k,d}^2 = \lim_{\varrho \to \infty} \frac{\text{Var}(Z_k^{(d)}(W_\varrho))}{|W_\varrho|} \qquad \text{for } k = 1, \ldots, d-1, \tag{5.2}$$

the existence of which is shown in the subsequent lemma.

**Lemma 5.1.** *Let there be given an $X/X_0$-nesting in $\mathbb{R}^d$ with stationary (not necessarily ergodic) initial tessellation $X$ and motion-invariant component tessellation $X_0$ satisfying the assumptions of Lemma 3.1 such that the asymptotic variance $\tau_{k,d}^2$ in the first part of Condition G exists. Then, the asymptotic variance $\widetilde{\tau}_{k,d}^2$ in (5.2) exists and takes the form*

$$\widetilde{\tau}_{k,d}^2 = (\tau_0^{(k,d)})^2 + (\lambda_0^{(kd)} m_k^{(d)})^2 \tau_{k,d}^2 \qquad \text{for } k = 1, \ldots, d-1, \tag{5.3}$$

*where*

$$(\tau_0^{(k,d)})^2 = \gamma \int_{\mathcal{P}_d^0} \text{Var}(\nu_{k-1}(X_0^{(d-1)} \cap C^{(k)})) P^0(\mathrm{d}C). \tag{5.4}$$

**Proof.** The proof of (5.3) is based on the representation of the variance of $Z_k^{(d)}(W_\varrho)$ for $k = 1, \ldots, d-1$ given in Lemma 3.1. From (3.6) and the first part of Condition G, it is easily seen that (5.3) holds if and only if the limit

$$\lim_{\varrho \to \infty} \frac{\mathbb{E}(\text{Var}_X(Z_k^{(d)}(W_\varrho)))}{|W_\varrho|}$$

$$= \lim_{\varrho \to \infty} \frac{\gamma}{|W_\varrho|} \int_{\mathcal{P}_d^0} \int_{\mathbb{R}^d} \text{Var}(\nu_{k-1}(X_0^{(d-1)} \cap C^{(k)} \cap (W_\varrho - x))) \, \mathrm{d}x P^0(\mathrm{d}C)$$

exists and equals $(\tau_0^{(k,d)})^2$, as defined in (5.4). To show this, we apply the same arguments as those already used in the proof of Lemma 3.1 to derive the estimate

$$\mathbb{E}(\text{Var}_X Z_k^{(d)}(W)) \leq \gamma |W| \int_{\mathcal{P}_d^0} \mathbb{E}\nu_{k-1}^2(X_0^{(d-1)} \cap C^{(k)}) P^0(\mathrm{d}C).$$

By multiple application of Fubini's theorem, we arrive at

$$\int_{\mathcal{P}_d^0} \int_{\mathbb{R}^d} \frac{\mathbb{E}\nu_{k-1}^2(X_0^{(d-1)} \cap C^{(k)} \cap (W_\varrho - x))}{|W_\varrho|} \, \mathrm{d}x P^0(\mathrm{d}C)$$

$$= \int_{\mathcal{P}_d^0} \mathbb{E}\Big( \int_{\mathbb{R}^d} \frac{\nu_{k-1}^2(X_0^{(d-1)} \cap C^{(k)} \cap (W_\varrho - x))}{|W_\varrho|} \, \mathrm{d}x \Big) P^0(\mathrm{d}C)$$

$$= \int_{\mathcal{P}_d^0} \mathbb{E}\Big( \int_{X_0^{(d-1)} \cap C^{(k)}} \int_{X_0^{(d-1)} \cap C^{(k)}} \frac{|(W_\varrho - u) \cap (W_\varrho - v)|}{|W_\varrho|} \nu_{k-1}(\mathrm{d}u)\nu_{k-1}(\mathrm{d}v) \Big) P^0(\mathrm{d}C)$$

$$\underset{\varrho \to \infty}{\longrightarrow} \int_{\mathcal{P}_d^0} \mathbb{E}\nu_{k-1}^2(X_0^{(d-1)} \cap C^{(k)}) P^0(\mathrm{d}C).$$

Note that, in view of $\lim_{\varrho \to \infty} |(W_\varrho - u) \cap (W_\varrho - v)|/|W_\varrho| = 1$ and $|(W_\varrho - u) \cap (W_\varrho - v)| \leq |W_\varrho|$, together with (3.5), we may apply Lebesgue's dominated convergence theorem. Likewise, we obtain that

$$\int_{\mathcal{P}_d^0} \int_{\mathbb{R}^d} \frac{(\mathbb{E}\nu_{k-1}(X_0^{(d-1)} \cap C^{(k)} \cap (W_\varrho - x)))^2}{|W_\varrho|} \, \mathrm{d}x P^0(\mathrm{d}C)$$

$$= (\lambda_0^{(k,d)})^2 \int_{\mathcal{P}_d^0} \int_{\mathbb{R}^d} \frac{\nu_k^2(C^{(k)} \cap (W_\varrho - x))}{|W_\varrho|} \, \mathrm{d}x P^0(\mathrm{d}C)$$

$$= (\lambda_0^{(k,d)})^2 \int_{\mathcal{P}_d^0} \int_{C^{(k)}} \int_{C^{(k)}} \frac{|(W_\varrho - u) \cap (W_\varrho - v)|}{|W_\varrho|} \nu_k(\mathrm{d}u)\nu_k(\mathrm{d}v) P^0(\mathrm{d}C)$$

$$\underset{\varrho \to \infty}{\longrightarrow} (\lambda_0^{(k,d)})^2 \int_{\mathcal{P}_d^0} \nu_k^2(C^{(k)}) P^0(\mathrm{d}C) = \int_{\mathcal{P}_d^0} (\mathbb{E}\nu_{k-1}(X_0^{(d-1)} \cap C^{(k)}))^2 P^0(\mathrm{d}C),$$

which completes the proof of Lemma 5.1. □

Under Condition F, we may decompose the normalized cumulative functionals given in (5.1), in analogy to (4.3), as

$$\frac{Z_k^{(d)}(W_\varrho) - \eta_k^{(d)}|W_\varrho|}{|W_\varrho|^{1/2}} = U_\varrho^{(k)} + m_k^{(d)} \mu_k^{(d)} V_\varrho^{(k)}, \tag{5.5}$$

where

$$U_\varrho^{(k)} = \frac{Z_k^{(d)}(W_\varrho) - \mathbb{E}_X(Z_k^{(d)}(W_\varrho))}{|W_\varrho|^{1/2}} \quad \text{and} \quad V_\varrho^{(k)} = \frac{\nu_k(X^{(k)} \cap W_\varrho) - \lambda_0^{(k,d)}|W_\varrho|}{|W_\varrho|^{1/2}}.$$

Notice the fact that, for $k = 1, \ldots, d-1$ and any fixed $\varrho > 0$, the random variables $U_\varrho^{(k)}$ and $V_\varrho^{(k)}$ are uncorrelated. Obviously, $V_\varrho^{(k)}$ is a (measurable) function of $X$ and

$\mathbb{E}_X(U_\varrho^{(k)}) = 0$ ($\mathbb{P}$-a.s), so

$$\mathbb{E}(U_\varrho^{(k)} V_\varrho^{(k)}) = \mathbb{E}(\mathbb{E}_X(U_\varrho^{(k)} V_\varrho^{(k)})) = \mathbb{E}(\mathbb{E}_X(U_\varrho^{(k)}) V_\varrho^{(k)}) = 0.$$

The following theorem states that, for any $k = 1, \ldots, d-1$, the two-dimensional vector $(U_\varrho^{(k)}, V_\varrho^{(k)})^\top$ converges in distribution to a mean zero Gaussian vector with independent components as $\varrho \to \infty$. This, in turn, implies the desired asymptotic normality of (5.1).

**Theorem 5.1.** *Consider an $X/X_0$-nesting in $\mathbb{R}^d$ observed through the increasing family of windows $W_\varrho$ with motion invariant component tessellation $X_0$ and stationary ergodic initial tessellation $X = \{\Xi_n\}_{n \geq 1}$ satisfying $\mathbb{E}D^d(\Xi^*) < \infty$, as well as (3.5), Conditions* F *and* G *Then,*

$$\begin{pmatrix} U_\varrho^{(k)} \\ V_\varrho^{(k)} \end{pmatrix} \xrightarrow[n \to \infty]{\mathrm{d}} \mathcal{N}\left( \begin{pmatrix} 0 \\ 0 \end{pmatrix}, \begin{pmatrix} (\tau_0^{(k,d)})^2 & 0 \\ 0 & \tau_{k,d}^2 \end{pmatrix} \right) \qquad \text{for } k = 1, \ldots, d-1. \tag{5.6}$$

*In particular, this implies that*

$$\frac{Z_k^{(d)}(W_\varrho) - \eta_k^{(d)}|W_\varrho|}{|W_\varrho|^{1/2}} \xrightarrow[n \to \infty]{} \mathcal{N}(0, \widetilde{\tau}_{k,d}^2) \qquad \text{for } k = 1, \ldots, d-1, \tag{5.7}$$

*where*

$$\eta_k^{(d)} = \lambda_0^{(k,d)} m_k^{(d)} \mu_k^{(d)} \quad \text{and} \quad \widetilde{\tau}_{k,d}^2 = (\tau_0^{(k,d)})^2 + (\lambda_0^{(k,d)} m_k^{(d)})^2 \tau_{k,d}^2.$$

**Proof.** We employ the method of characteristic functions (cf., e.g., [17] for details). Hence, we must show that the characteristic function $f_\varrho(s,t)$ of the random vector $(U_\varrho^{(k)}, V_\varrho^{(k)})^\top$ defined by

$$f_\varrho(s,t) = \mathbb{E} \exp\{\mathrm{i}sU_\varrho^{(k)} + \mathrm{i}tV_\varrho^{(k)}\}$$

converges to the characteristic function of the Gaussian random vector that occurs as limit in (5.6), that is,

$$f_\varrho(s,t) \xrightarrow[\varrho \to \infty]{} \exp\left\{ -\frac{s^2}{2}(\tau_0^{(k,d)})^2 - \frac{t^2}{2}\tau_{k,d}^2 \right\} \qquad \text{for all } s, t \in \mathbb{R}^1.$$

For this, we introduce the decomposition $f_\varrho(s,t) = \sum_{i=1}^3 f_\varrho^{(i)}(s,t)$, where

$$f_\varrho^{(1)}(s,t) = \mathbb{E}\left[ \mathbb{E}_X\left( \exp\{\mathrm{i}sU_\varrho^{(k)}\} - \exp\left\{ -\frac{s^2}{2|W_\varrho|}\mathrm{Var}_X(Z_k^{(d)}(W_\varrho)) \right\} \right) \exp\{\mathrm{i}tV_\varrho^{(k)}\} \right],$$

$$f_\varrho^{(2)}(s,t) = \mathbb{E}\left[ \left( \exp\left\{ -\frac{s^2}{2|W_\varrho|}\mathrm{Var}_X(Z_k^{(d)}(W_\varrho)) \right\} - \exp\left\{ -\frac{s^2}{2}(\tau_0^{(k,d)})^2 \right\} \right) \exp\{\mathrm{i}tV_\varrho^{(k)}\} \right]$$

and

$$f_\varrho^{(3)}(s,t) = \exp\left\{-\frac{s^2}{2}(\tau_0^{(k,d)})^2\right\} \mathbb{E}\exp\{\mathrm{it}V_\varrho^{(k)}\}.$$

In view of Condition G, the continuity theorem for (one-dimensional) characteristic functions yields that

$$\lim_{\varrho\to\infty} f_\varrho^{(3)}(s,t) = \exp\left\{-\frac{s^2}{2}(\tau_0^{(k,d)})^2\right\} \lim_{\varrho\to\infty} \mathbb{E}\exp\{\mathrm{it}V_\varrho^{(k)}\}$$

$$= \exp\left\{-\frac{s^2}{2}(\tau_0^{(k,d)})^2 - \frac{t^2}{2}\tau_{k,d}^2\right\}$$

for all $s,t \in \mathbb{R}^1$. Hence, it remains to prove $f_\varrho^{(i)}(s,t) \longrightarrow_{\varrho\to\infty} 0$ for $i = 1,2$. For this, we subsequently show that

$$\frac{\mathrm{Var}_X(Z_k^{(d)}(W_\varrho))}{|W_\varrho|} = \frac{1}{|W_\varrho|}\sum_{n\geq 1}\mathrm{Var}_X(\nu_{k-1}(X_0^{(d-1)}\cap\Xi_n^{(k)}\cap W_\varrho)) \xrightarrow[\varrho\to\infty]{\mathrm{a.s.}} (\tau_0^{(k,d)})^2, \quad (5.8)$$

where the first equality follows from (3.9) and the almost sure convergence of the argued sum in (5.8) can be argued by some modified ergodic theorem for random tessellations (cf. Theorem 4.1 in [14]), which states that

$$\frac{1}{|W_\varrho|}\sum_{n\geq 1}\mathbb{1}_{\{\Xi_n\cap W_\varrho\neq\emptyset\}}g(\Xi_n\cap W_\varrho) \xrightarrow[\varrho\to\infty]{\mathrm{a.s.}} \gamma\mathbb{E}g(\Xi^*) \quad (5.9)$$

for any $\mathcal{B}(\mathcal{F}_d')$-measurable, translation invariant set function $g(\cdot)$ defined on sets of the form $C\cap W$, where $C$ is a $d$-polytope and $W\in\mathcal{C}_d'$, and satisfying the monotonicity property $g(C\cap W)\leq g(C\cap W')$ for $W\subseteq W'$. It is easily verified, by checking the proof of Theorem 4.1 in [14], that these restrictions imposed on $g(\cdot)$, together with $\mathbb{E}D^d(\Xi^*) < \infty$ and $\mathbb{E}g(\Xi^*) < \infty$, suffice for (5.9) to hold. Applying (5.9) to $g_1(C\cap W) = \mathbb{E}\nu_{k-1}^2(X_0^{(d-1)}\cap C^{(k)}\cap W)$ and $g_2(C\cap W) = (\lambda_0^{(k,d)})^2\nu_k^2(C^{(k)}\cap W)$, where $C^{(k)}$ denotes the union of $k$-faces of the $d$-polytope $C$, we see that (5.9) also holds for $g = g_1 - g_2$. Thus, (5.8) is proved. From (5.8) and Lebesgue's dominated convergence theorem, we conclude that

$$|f_\varrho^{(2)}(s,t)| \leq \mathbb{E}\left|\exp\left\{-\frac{s^2}{2|W_\varrho|}\mathrm{Var}_X(Z_k^{(d)}(W_\varrho))\right\} - \exp\left\{-\frac{s^2}{2}((\tau_0)^{(k,d)})^2\right\}\right| \xrightarrow[\varrho\to\infty]{} 0. \quad (5.10)$$

To show that $f_\varrho^{(1)}(s,t)$ becomes arbitrarily small as $W_\varrho$ grows large, we start with the obvious estimate

$$|f_\varrho^{(1)}(s,t)| \leq \mathbb{E}\left|\mathbb{E}_X\exp\{\mathrm{is}U_\varrho^{(k)}\} - \exp\left\{-\frac{s^2}{2|W_\varrho|}\mathrm{Var}_X(Z_k^{(d)}(W_\varrho))\right\}\right| \quad (5.11)$$

for $k = 1, \ldots, d-1$. Next, we express $U_\varrho^{(k)}$ in terms of the centered measures $\theta_n^{(k)}(\cdot) = \nu_{k-1}(X_n^{(d-1)} \cap \Xi_n^{(k)} \cap (\cdot)) - \lambda_0^{(k,d)}\nu_k(\Xi_n^{(k)} \cap (\cdot))$, which are conditionally independent given the initial tessellation $X = \{\Xi_n\}_{n \geq 1}$. We have

$$U_\varrho^{(k)} = \frac{1}{|W_\varrho|^{1/2}} \sum_{n \geq 1} \mathbb{1}_{\{W_\varrho \cap \Xi_n \neq \emptyset\}} \theta_n^{(k)}(W_\varrho)$$

and introduce, for fixed $\delta \geq 0$, the conditional Lindeberg function

$$L_\varrho^{(k)}(\delta) = \frac{1}{|W_\varrho|} \sum_{n \geq 1} \mathbb{1}_{\{W_\varrho \cap \Xi_n \neq \emptyset\}} \mathbb{E}_X(\theta_n^{(k)}(W_\varrho))^2 \mathbb{1}_{\{|\theta_n^{(k)}(W_\varrho)| \geq \delta |W_\varrho|^{1/2}\}}. \qquad (5.12)$$

Further, for any $\varepsilon > 0$ and $\delta > 0$, we define the events

$$G_\varrho(\mathrm{e}, \delta) = \{L_\varrho^{(k)}(\delta) \leq \varepsilon\} \quad \text{and} \quad H_\varrho(\varepsilon) = \{|L_\varrho^{(k)}(0) - (\tau_0^{(k,d)})^2| \leq \varepsilon\}.$$

Since $L_\varrho^{(k)}(0) = \mathrm{Var}_X(Z_k^{(d)}(W_\varrho))/|W_\varrho|$, it follows from (5.8) that $\mathbb{P}(H_\varrho^c(\varepsilon)) \longrightarrow_{\varrho \to \infty} 0$. Below, we also need that $\mathbb{P}(G_\varrho^c(\varepsilon, \delta)) \longrightarrow_{\varrho \to \infty} 0$, following from the stronger result $L_\varrho(\delta) \xrightarrow{\mathrm{a.s.}}_{\varrho \to \infty} 0$, which we can show in the following way. Replacing $\theta_n^{(k)}(W_\varrho)$ in (5.12) by $\xi_n^{(k)}(\Xi_n^{(k)} \cap W_\varrho) = \nu_{k-1}(X_n^{(d-1)} \cap \Xi_n^{(k)} \cap W_\varrho) + \lambda_0^{(k,d)} \nu_k(\Xi_n^{(k)} \cap W_\varrho)$ leads to the inequality

$$L_\varrho^{(k)}(a|W_\varrho|^{-1/2}) \leq \frac{1}{|W_\varrho|} \sum_{n \geq 1} \mathbb{1}_{\{W_\varrho \cap \Xi_n \neq \emptyset\}} \mathbb{E}_X(\xi_0^{(k)}(\Xi_n^{(k)} \cap W_\varrho))^2 \mathbb{1}_{\{\xi_0^{(k)}(\Xi_n^{(k)} \cap W_\varrho) \geq a\}}$$

for any $a > 0$. The set function $g(C \cap W) = \mathbb{E}_X(\xi_0^{(k)}(C^{(k)} \cap W))^2 \mathbb{1}_{\{\xi_0^{(k)}(C^{(k)} \cap W) \geq a\}}$ is translation invariant (due to the stationarity of $X_0$) and increases whenever $W$ expands. Hence, $g(\cdot)$ fulfills the conditions needed to establish the almost sure convergence in (5.9). This implies that

$$\mathbb{P}\left(\limsup_{\varrho \to \infty} L_\varrho(a|W_\varrho|^{-1/2}) \leq \gamma \int_{\mathcal{P}_d^0} \mathbb{E}_X(\xi_0^{(k)}(C^{(k)}))^2 \mathbb{1}_{\{\xi_0^{(k)}(C^{(k)}) \geq a\}} P^0(\mathrm{d}C)\right) = 1.$$

Consequently, by (3.5), $\mathbb{P}(\limsup_{\varrho \to \infty} L_\varrho^{(k)}(\delta) \leq \varepsilon) = 1$ for any $\varepsilon > 0$, that is, $L_\varrho^{(k)}(\delta) \xrightarrow{\mathrm{a.s.}}_{\varrho \to \infty} 0$. A suitable upper bound of the right-hand side of (5.11) can be obtained when both events $G_\varrho(\varepsilon, \delta)$ and $H_\varrho(\varepsilon)$ occur. From (5.11), it is easily seen that $|f_\varrho^{(1)}(s,t)|$ does not exceed the sum

$$\mathbb{E}\mathbb{1}_{\{G_\varrho(\varepsilon, \delta) \cap H_\varrho(\varepsilon)\}} \left| \mathbb{E}_X \exp\{\mathrm{i}sU_\varrho^{(k)}\} - \exp\left\{-\frac{s^2}{2} L_\varrho(0)\right\} \right| + 2\mathbb{P}(G_\varrho^c(\varepsilon, \delta) \cup H_\varrho^c(\varepsilon)). \quad (5.13)$$

We proceed with the factorization of the conditional characteristic function of $U_\varrho^{(k)}$ given $X$, using the conditional independence of the random variables $\theta_n^{(k)}(W_\varrho)$, and

obtain

$$\mathbb{E}_X \exp\{\mathrm{i}sU_\varrho^{(k)}\} = \prod_{n\geq 1} \mathbb{E}_X \exp\{\mathrm{i}s|W_\varrho|^{-1/2}\theta_n^{(k)}(W_\varrho)\}.$$

Expressing the first equality in (5.8) by the centered measures $\theta_n^{(k)}(\cdot)$, we have

$$\exp\left\{-\frac{s^2}{2|W_\varrho|}\operatorname{Var}_X(Z_k^{(d)}(W_\varrho))\right\} = \prod_{n\geq 1}\exp\left\{-\frac{s^2}{2|W_\varrho|}\mathbb{E}_X(\theta_n^{(k)}(W_\varrho))^2\right\}.$$

By means of the elementary inequality $|x_1\cdots x_n - y_1\cdots y_n| \leq |x_1 - y_1| + \cdots + |x_n - y_n|$ for complex numbers $x_i, y_i$ lying on the unit disc, we arrive at the estimate

$$\begin{aligned}
&\left|\mathbb{E}_X\exp\{\mathrm{i}sU_\varrho^{(k)}\} - \exp\left\{-\frac{s^2}{2|W_\varrho|}\operatorname{Var}_X(Z_k^{(d)}(W_\varrho))\right\}\right| \\
&\leq \sum_{n\geq 1}\mathbb{1}_{\{\Xi_n\cap W_\varrho\neq\emptyset\}}\left|\mathbb{E}_X\exp\left\{\mathrm{i}s\frac{\theta_n^{(k)}(W_\varrho)}{\sqrt{|W_\varrho|}}\right\} - \exp\left\{-\frac{s^2}{2|W_\varrho|}\mathbb{E}_X(\theta_n^{(k)}(W_\varrho))^2\right\}\right|.
\end{aligned} \tag{5.14}$$

Further, using the well-known inequality $|e^{\mathrm{i}x} - \sum_{k=0}^{n-1}\frac{(\mathrm{i}x)^k}{k!}| \leq \frac{|x|^n}{n!}$ (with $x\in\mathbb{R}^1$) for $n=2$ and $n=3$, we find that, for any $\delta > 0$,

$$\begin{aligned}
&\left|\mathbb{E}_X\left(\exp\left\{\mathrm{i}s\frac{\theta_n^{(k)}(W_\varrho)}{\sqrt{|W_\varrho|}}\right\} - 1 - \mathrm{i}s\frac{\theta_n^{(k)}(W_\varrho)}{\sqrt{|W_\varrho|}} + \frac{s^2}{2|W_\varrho|}\mathbb{E}_X(\theta_n^{(k)}(W_\varrho))^2\right)\right| \\
&\leq \frac{s^2}{|W_\varrho|}\mathbb{E}_X(\theta_n^{(k)}(W_\varrho))^2\mathbb{1}_{\{|\theta_n^{(k)}(W_\varrho)|\geq\delta\sqrt{|W_\varrho|}\}} \\
&\quad + \frac{|s|^3}{6|W_\varrho|^{3/2}}\mathbb{E}_X|\theta_n^{(k)}(W_\varrho)|^3\mathbb{1}_{\{|\theta_n^{(k)}(W_\varrho)|\leq\delta\sqrt{|W_\varrho|}\}}.
\end{aligned} \tag{5.15}$$

Analogously, applying the inequality $|e^{-x} - 1 + x| \leq x^2/2$ for $x\geq 0$ gives

$$\begin{aligned}
&\left|\exp\left\{-\frac{s^2}{2|W_\varrho|}\mathbb{E}_X(\theta_n^{(k)}(W_\varrho))^2\right\} - 1 + \frac{s^2}{2|W_\varrho|}\mathbb{E}_X(\theta_n^{(k)}(W_\varrho))^2\right| \\
&\leq \frac{s^4}{4|W_\varrho|^2}(\mathbb{E}_X(\theta_n^{(k)}(W_\varrho))^2)^2 \leq \frac{s^4}{4|W_\varrho|}\mathbb{E}_X(\theta_n^{(k)}(W_\varrho))^2(\delta^2 + L_\varrho^{(k)}(\delta)),
\end{aligned} \tag{5.16}$$

where we have used, in addition, that, for any $n\geq 1$ and $\delta > 0$,

$$\mathbb{E}_X(\theta_n^{(k)}(W_\varrho))^2 \leq \delta^2|W_\varrho| + \mathbb{E}_X(\theta_n^{(k)}(W_\varrho))^2\mathbb{1}_{\{|\theta_n^{(k)}(W_\varrho)|\geq\delta\sqrt{|W_\varrho|}\}} \leq |W_\varrho|(\delta^2 + L_\varrho^{(k)}(\delta)).$$

Finally, combining the above estimates (5.14), (5.15) and (5.16) and taking into account both the abbreviation (5.12) and the fact that $\mathbb{E}_X \theta_n^{(k)}(W_\varrho) = 0$, we find that

$$\left| \mathbb{E}_X \exp\{isU_\varrho^{(k)}\} - \exp\left\{-\frac{s^2}{2} L_\varrho^{(k)}(0)\right\} \right| \leq s^2 L_\varrho^{(k)}(\delta) + \frac{|s|^3 \delta}{6} L_\varrho^{(k)}(0) + \frac{s^4}{2} L_\varrho^{(k)}(0)(\delta^2 + L_\varrho^{(k)}(\delta)).$$

Regarding the latter inequality on the event $\{G_\varrho(\varepsilon, \delta) \cap H_\varrho(\varepsilon)\}$, we obtain from (5.13) that

$$\limsup_{\varrho \to \infty} |f_\varrho^{(1)}(s,t)| \leq s^2 \varepsilon + \left(\frac{|s|^3 \delta}{6} + \frac{s^4}{2}(\delta^2 + \varepsilon)\right)(\varepsilon + (\tau_0^{(k,d)})^2)$$

for arbitrary $\varepsilon, \delta > 0$. Thus, $\lim_{\varrho \to \infty} f_\varrho^{(1)}(s,t) = 0$ which completes the proof of (5.6). The proof of Theorem 5.1 ends with an application of (5.6) and the continuous mapping theorem (cf. [17]) to the linear combination (5.5), which proves (5.7). □

## 6. Examples of weakly dependent random tessellations

There are only a few papers (e.g., [3, 10]) concerning weak dependence properties of stationary random tessellations apart from ergodicity. In fact, the assumption of ergodicity turns out to be the weakest form of asymptotic independence of distant parts of a stationary tessellation $X = \{\Xi\}_{n \geq 1}$. Due to the individual spatial ergodic theorem (cf. [6]) ergodicity guarantees strong consistency for a number of intensity estimators based on a single observation in an expanding sampling window. To establish asymptotic normality of these estimators, the distribution of $X$ must satisfy certain *mixing conditions* expressed in terms of corresponding *mixing coefficients*.

In the context of random tessellations $X = \{\Xi\}_{n \geq 1}$ in $\mathbb{R}^d$, the $\alpha$- and $\beta$-mixing condition have proved meaningful with mixing coefficients defined by

$$\begin{aligned}
\alpha(\mathcal{A}_X(F_1), \mathcal{A}_X(F_2)) &= \sup_{A_1 \in \mathcal{A}_X(F_1), A_2 \in \mathcal{A}_X(F_2)} |\mathbb{P}(A_1 \cap A_2) - \mathbb{P}(A_1)\mathbb{P}(A_2)|, \\
\beta(\mathcal{A}_X(F_1), \mathcal{A}_X(F_2)) &= \mathbb{E} \sup_{A_2 \in \mathcal{A}_X(F_2)} |\mathbb{P}(A_2|\mathcal{A}_X(F_1)) - \mathbb{P}(A_2)|,
\end{aligned}$$
(6.1)

where $F_1, F_2$ are disjoint closed subsets of $\mathbb{R}^d$ and $\mathcal{A}_X(F)$ denotes the $\sigma$-algebra generated by the random closed set $(\bigcup_{n \geq 1} \partial \Xi_n \cap F$ in the sense of Matheron (cf. [20] and also [10]). It is easily verified that $\alpha(\mathcal{A}_X(F_1), \mathcal{A}_X(F_2)) \leq \beta(\mathcal{A}_X(F_1), \mathcal{A}_X(F_2))$. However, the behaviour of both mixing coefficients is nearly the same for most of the models in stochastic geometry when the distance between $F_1$ and $F_2$ becomes large (cf. [12]). To verify Condition G, we are faced with two problems. First, to find, from the model assumptions, sharp bounds on the above mixing coefficients for $F_1 = F(a) := [-a,a]^d$ and $F_2 := G(b) = \mathbb{R}^d \setminus (-b,b)^d$ for $b > a$ and, second, to prove a suitable CLT (or to use a known CLT) for weakly dependent random fields whose assumptions follow from the derived estimates. For more details on CLTs of random fields and mixing conditions, including the influence of the dimension, the reader is referred to [23].

In [10], the $\beta$-mixing coefficient $\beta(\mathcal{A}_X(F(a)), \mathcal{A}_X(G(a+r)))$ could be estimated for Voronoi tessellations in terms of the $\beta$-mixing coefficient and certain void probabilities of the generating stationary point process of nuclei. In the special case of Voronoi tessellations generated by Poisson cluster processes with cluster radius $R_0$ satisfying $\mathbb{E} \exp\{h R_0\} < \infty$ for some $h > 0$, the general bound decays exponentially in $r$. More precisely, it can be shown that

$$\beta(\mathcal{A}_X(F(a)), \mathcal{A}_X(G(a+r))) \le c_1 \left( \left( \frac{r}{a} \right)^{d-1} + \left( \frac{a}{r} \right)^{d-1} \right) \exp\{-c_2 r\} \quad (6.2)$$

for any $r \ge 1$ and $a \ge 1/2$, where the positive constants $c_1, c_2$ depend only on the dimension $d$, $h > 0$ and the intensity of the Poisson process of cluster centers. An estimate similar to (6.2) holds for Poisson soft-core processes (cf. [29]), provided the soft-core radius $R_0$ possesses an exponential moment. Furthermore, Condition G could be verified in [10] for stationary random tessellations $X = \{\Xi\}_{n \ge 1}$ in $\mathbb{R}^d$ satisfying $\mathbb{E}(\nu_k(X^{(k)} \cap [0,1)^d))^{2+\delta} < \infty$ for $1 \le k \le d-1$ and some $\delta > 0$ and

$$\beta(\mathcal{A}_X(F(a)), \mathcal{A}_X(G(a+r))) \le a^{d-1} \beta_1(r) \mathbb{1}_{[1,ca]}(r) + \beta_2(r) \mathbb{1}_{(ca,\infty)}(r)$$

for any $a \ge 1/2$ and $r \ge 1$, where $c \ge 2$ is a constant independent of both $a$ and $r$. Furthermore, $\beta_1(\cdot)$ and $\beta_2(\cdot)$ are non-increasing functions on $[1,\infty)$ such that

$$r^{2d-1} \beta_1(r) \underset{\varrho \to \infty}{\longrightarrow} \quad \text{and} \quad \sum_{r \ge 1} r^{d-1} (\beta_2(r))^{\delta/(2+\delta)} < \infty.$$

Hence, from (6.2), we obtain, for a Voronoi tessellation $X = \{\Xi\}_{n \ge 1}$ in $\mathbb{R}^d$ generated by a stationary Poisson process with intensity $\gamma > 0$, that

$$|W_\varrho|^{-1/2} (\nu_k(X^{(k)} \cap W_\varrho) - \overline{\mu}_k^{(d)}(\nu) \gamma^{(d-k)/d} |W_\varrho|) \underset{n \to \infty}{\overset{\mathrm{d}}{\longrightarrow}} \mathcal{N}(0, \overline{\tau}_{k,d}^2(\nu) \gamma^{(d-2k)/d}) \quad (6.3)$$

for $k = 1, \ldots, d-1$, where the mean value $\overline{\mu}_k^{(d)}(\nu) = \mathbb{E}(\nu_k(X^{(k)} \cap [0,1)^d)$ and the asymptotic variance $\overline{\tau}_{k,d}^2(\nu)$ refer to a PVT with intensity $\gamma = 1$. The scaling rates in (6.3) are easily seen from the scaling property of the stationary Poisson process in $\mathbb{R}^d$ giving $\nu_k(X^{(k)} \cap W_\varrho) = \gamma^{-k/d} \nu_k(\overline{X}^{(k)} \cap W_{\varrho \gamma^{1/d}})$, where $\overline{X}^{(k)}$ denotes the union of $k$-facets of a PVT with unit intensity. An explict formula for the intensity $\overline{\mu}_k^{(d)}(\nu)$ was first found (by R. Miles) to be

$$\overline{\mu}_k^{(d)}(\nu) = \frac{(2\pi)^{d-k+1} \Gamma(d-k+k/d)}{(d-k+1)! d} \frac{\kappa_{d(d-k)+k-2}}{\kappa_{d(d-k)+k-1}} \frac{\kappa_{k-1}}{(\kappa_d)^{k/d}} \left( \frac{\kappa_{d-1}}{\kappa_d} \right)^{d-k}$$

and derived in a different way in [21], page 64, whereas, for $\overline{\tau}_{k,d}^2(\nu)$, no analytic expression is currently known. In the planar case, it is well known that $\overline{\mu}_1^{(2)}(\nu) = 2$ and the approximate value $\overline{\tau}_{1,2}^2(\nu) = 1.0445685$ was found in [1] by numerical evaluation of rather involved multiple integrals.

We are now in a position to establish the CLT in (5.7) in a more explict form for the case of a planar $X/X_0$-nesting with initial tessellation $X$ being a PVT with cell intensity $\gamma > 0$ and component tessellation $X_0$ being either a PLT generated by a motion-invariant Poisson line process with intensity $\lambda > 0$ or another PVT with cell intensity $\lambda > 0$. In both cases, we have

$$\frac{Z_1^{(2)}(W_\varrho) - \eta_1^{(2)}|W_\varrho|}{|W_\varrho|^{1/2}} \xrightarrow[n \to \infty]{\mathrm{d}} \mathcal{N}(0, \widetilde{\tau}_{1,2}^2), \tag{6.4}$$

where

$$\eta_1^{(2)} = 4\lambda_0^{(1,2)}\sqrt{\lambda} \quad \text{and} \quad \widetilde{\tau}_{1,2}^2 = (\tau_0^{(1,2)})^2 + 1.6934(\lambda_0^{(1,2)})^2.$$

From (3.2) we get $\lambda_0^{(1,2)} = c_1^{(2)}\lambda_0^{(2,2)}$ with $c_1^{(2)} = 2/\pi$ and $\lambda_0^{(2,2)} = 2\sqrt{\lambda}$ if $X_0$ is a PVT (cf. [24], page 314 and $\lambda_0^{(2,2)} = \lambda$ if $X_0$ is a PLT (cf (2.8)). To calculate

$$(\tau_0^{(1,2)})^2 = \gamma \int_{\mathcal{P}_d^0} \mathrm{Var}(\nu_0(X_0^{(1)} \cap C^{(1)})) P^0(\mathrm{d}C),$$

we first consider the case when $X_0$ is a PLT with intensity $\lambda$. Then, $\nu_0(X_0^{(1)} \cap C^{(1)})$ equals twice the number $N(C)$ of Poisson lines hitting the polygon $C$. It is well known that $N(C)$ is Poisson distributed with mean (and variance) $\lambda P(C)/\pi$, where $P(C)$ denotes the perimeter of $C$. Hence, by $\mathbb{E}P(\Xi^*) = 4\gamma^{-1/2}$ (cf. [24], page 314) we obtain

$$(\tau_0^{(1,2)})^2 = \gamma \frac{4\lambda \mathbb{E}P(\Xi^*)}{\pi} = \frac{16}{\pi}\sqrt{\gamma}\lambda.$$

If $X_0$ is a PVT with cell intensity $\lambda$, we may again exploit the scaling properties of PVTs, giving

$$(\tau_0^{(1,2)})^2 = \sqrt{\gamma\lambda}\mathbb{E}\,\mathrm{Var}_{\overline{X}}(\nu_0(\overline{X}_0^{(1)} \cap \overline{\Xi}^{*,(1)})),$$

where $\overline{X}$ and $\overline{X}_0$ are independent planar PVTs, both with unit cell intensity, and where $\overline{X}_0^{(1)} \cap \overline{\Xi}^{*,(1)}$ denotes the finite set of points on the boundary of the typical cell of $\overline{X}$ induced by the 1-facets of $\overline{X}_0$. A large-scale simulation study yields $\mathbb{E}\,\mathrm{Var}_{\overline{X}}(\nu_0(\overline{X}_0^{(1)} \cap \overline{\Xi}^{*,(1)})) = 2.7023$. Summarizing the above results, we obtain the following expressions for $\eta_1^{(2)}$ and $\widetilde{\tau}_{1,2}^2$ in (6.4), namely,

$$\eta_1^{(2)} = \frac{8}{\pi}\sqrt{\gamma}\lambda \quad \text{and} \quad \widetilde{\tau}_{1,2}^2 = \frac{16}{\pi}\sqrt{\gamma}\lambda + 1.6934\lambda^2 \qquad \text{(if } X_0 \text{ is a PLT with intensity } \lambda),$$

$$\eta_1^{(2)} = \frac{16}{\pi}\sqrt{\gamma\lambda} \quad \text{and} \quad \widetilde{\tau}_{1,2}^2 = 2.7023\sqrt{\gamma\lambda} + 6.7736\lambda \qquad \text{(if } X_0 \text{ is a PVT with intensity } \lambda).$$

To conclude, it should be mentioned that Condition G can also be verified for a large class of Laguerre tessellations generated by Poisson-based point processes. The values of the variances in the previous formulae for higher dimensions can only be obtained by extensive simulation studies. Several generalizations of Theorems 4.1 and 5.1 are possible. For example, the manifold process $X_0^{(d-1)}$ of $(d-1)$-facets can be replaced by the union of $k$-facets of $X_0$ for $1 \leq k \leq d-2$. In case of anisotropic component tessellations, the rose of directions of $X_0$ is needed to express the mean and variance of $Z_k^{(d)}(W_\varrho)$.

# Acknowledgements

The authors are grateful to the referees for their valuable comments and suggestions. This research was supported by France Telecom R&D through research Grant 42366897.